\documentclass[12pt]{amsart}
\usepackage[utf8]{inputenc}
\usepackage{amsmath}
\usepackage{amssymb}
\usepackage{xcolor}
\usepackage{fullpage}
\usepackage[colorlinks=true, pdfstartview=FitH, linkcolor=blue, citecolor=blue, urlcolor=blue]{hyperref}

\theoremstyle{plain}
\newtheorem{theorem}{Theorem}[section]
\newtheorem{lemma}[theorem]{Lemma}

\newtheorem{proposition}[theorem]{Proposition}

\theoremstyle{definition}
\newtheorem{definition}[theorem]{Definition}

\newtheorem{remark}[theorem]{Remark}

\newtheorem*{theoremA*}{Theorem A}

\newcommand{\cA}{{\mathcal{A}}}
\newcommand{\cB}{{\mathcal{B}}}

\DeclareMathOperator{\eps}{\varepsilon}
\renewcommand{\pmod}[1]{{\ifmmode\text{\rm\ (mod~$#1$)}\else\discretionary{}{}{\hbox{ }}\rm(mod~$#1$)\fi}}

\title{Explicit Burgess inequalities for cubefree moduli}

\begin{document}

\author[E. Hasanalizade]{Elchin Hasanalizade}
\address[Elchin Hasanalizade]{School of Information Technologies and Engineering\\ADA University, Baku, Azerbaijan}
\email{ehasanalizade@ada.edu.az}

\author[H. Lin]{Hua Lin}
\address[Hua Lin]{Department of Mathematics\\
Northwestern University, Evanston, IL, USA}
\email{hua.lin@northwestern.edu}

\author[A. Luna Mart\'inez]{Andradis Luna Mart\'inez}
\address[Andradis Luna Mart\'inez]{Escuela de Matem\'aticas\\Universidad Aut\'onoma de Santo Domingo, Santo Domingo, Dominican Republic}
\email{aluna71@uasd.edu.do}

\author[G. Martin]{Greg Martin}
\address[Greg Martin]{Department of Mathematics\\
University of British Columbia\\
Room 121, 1984 Mathematics Road\\
Vancouver, BC, Canada  V6T 1Z2}
\email{gerg@math.ubc.ca}

\author[E. Trevi{\~n}o]{Enrique Trevi{\~n}o}
\address[Enrique Trevi\~no]{Department of Mathematics and Computer Science\\Lake Forest College, Lake Forest, IL, USA}
\email{trevino@lakeforest.edu}

\begin{abstract}
    Burgess proved that for $\chi_q$ a primitive Dirichlet character modulo $q$ with $q$ cubefree, $\sum_{M< n\le M+N}\chi_q(n)= O\left(N^{1-\frac{1}{r}}q^{\frac{r+1}{4r^2}+\epsilon}\right)  $ for all integers $r\ge1.$ More recently, explicit versions with prime moduli $q$ were computed by Booker, McGown, Trevi\~{n}o, and Francis, with applications to finding the least $k$-th power residue, and bounding the size of Dirichlet $L$-functions just to name a few. Jain-Sharma, Khale, and Liu proved an explicit estimate for $r=2.$ We improve their explicit constant for $r = 2$ and compute an explicit Burgess bound for cubefree $q$ for $r\ge 3$.
\end{abstract}

\maketitle

\section{Introduction}
Let $q$ be a positive integer.
Given a primitive Dirichlet character $\chi\pmod q$ and integers~$M$ and~$N\ge1$, we consider the character sum
\begin{equation*}
    S_\chi\left(M, N\right)=\sum_{M<n \le M+N}\chi(n).
\end{equation*}

Burgess, in a series of papers~\cite{burgess:1962,Bur1962,Bur19632,Bur1963,burgess86}, proved
$$|S_{\chi}(M,N)|\le C_{\varepsilon,r} N^{1-\frac{1}{r}}q^{\frac{r+1}{4r^2}+\eps},$$
for some constant $C_{\varepsilon,r}$ that depends only on $\varepsilon$ and $r$ for $r = 2,3$ and for any $r\ge 4$ when $q$ is cubefree. In recent years, there has been progress in obtaining explicit estimates for $S_{\chi}(M,N)$. In the case when $q$ is prime, Booker \cite{Boo2006} proved a result for $\chi$ quadratic, and McGown \cite{mcgown2012} and Trevi\~no \cite{ET2} made improvements for general $\chi$. The current best result for $q$ prime and $2\le r\le 10$ is due to Francis \cite{francis2021}. An example of an explicit Burgess inequality for $q = p$ prime is the following uniform bound (from \cite{ET2})
$$|S_{\chi}(M,N)| \le 2.74 N^{1-\frac{1}{r}}p^{\frac{r+1}{4r^2}}(\log{p})^{\frac{1}{r}},$$
for all integers $r\ge 2$ and for all primes $p\ge 10^7$.

These explicit results have been used to determine the class number of a quadratic field of large discriminant \cite{Boo2006}, to bound the least quadratic non-residue modulo a prime \cite{francis2021,ET2}, to classify Norm-Euclidean cyclic fields \cite{mcgown2012}, to bound the least primitive root modulo a prime \cite{COT-2016}, to give an explicit bound for $L(1,\chi)$ in the case of quadratic characters \cite{JRT-2023}, to compute the rank of an elliptic curve with large rank \cite{Elkies}, among other applications.

In the case of composite moduli $q$, the only explicit result is the following due to Jain-Sharma, Khale, and Liu \cite{JKL}.
\begin{theoremA*} \label{TeoA}
    Let $\chi$ be a primitive character modulo $q$ with $q \ge e^{e^{9.594}}$. Then, for $N\le q^{\frac{5}{8}}$,
    $$\left|S_{\chi}(M,N)\right|\le 9.07 \sqrt{N} q^{\frac{3}{16}} \log^{\frac{1}{4}}(q)\left(2^{\omega(q)}\tau(q)\right)^{\frac{3}{4}}\left(\frac{q}{\phi(q)}\right)^{\frac{1}{2}},$$
    where $\tau(q)$ is the number of divisors of $q$ and $\omega(q)$ is the number of distinct prime factors of~$q$.
\end{theoremA*}

In other words, the only explicit result for composite moduli $q$ is for the case $r =2$. We improve this result for $r =2$ and $q$ large enough, and we also generalize the result to all $r\ge 3$ when $q$ is cubefree. 

To state our theorem, we use the following standard labels of arithmetic functions:
let $\omega(q)$ denote the number of distinct prime factors of~$q$, let $\phi(q)$ be Euler's totient function, and let $\tau_{k}(q)$ denote the number of ways to write $q$ as an ordered product of~$k$ integers, so that $\tau(n) = \tau_2(n)$ is the ordinary number-of-divisors function. Also define
\begin{equation} \label{meet Mr. Q}
m_r(q) = \min \biggl\{ \tau_{2r}(q), \biggl( \frac{\tau(q)}2 \biggr)^{2r-1} , \frac{q}{2r} \biggr\}.
\end{equation}

\begin{theorem}\label{main thm 1}
Let $r\ge 2$ be an integer and $\chi$ be a primitive Dirichlet character modulo~$q$. Let $C(r)$ be defined as in Table \ref{Table 1}. Let $a(r) = 2\log{2}\left(3.0758r+1.38402\log(4r)-1.5379\right)$. Then, for $q \ge \max\{10^{1143},e^{e^{a(r)}}\}$, if $r=2$ or $q$ is cubefree, we have
    $$|S_{\chi}(M,N)| \le C(r)N^{1-\frac{1}{r}} q^{\frac{r+1}{4r^2}}(\log{q})^{\frac{1}{2r}}\left((4r)^{\omega(q)}m_r(q)\right)^{\frac{1}{2r}-\frac{1}{2r^2}}\left(\frac{q}{\phi(q)}\right)^{\frac{1}{r}}.$$

    Furthermore, as $q\rightarrow\infty$, we have a constant $D(r)$ from Table \ref{Table 1} such that
    $$|S_{\chi}(M,N)| \le (D(r)+o(1))N^{1-\frac{1}{r}} q^{\frac{r+1}{4r^2}}(\log{q})^{\frac{1}{2r}}\left((4r)^{\omega(q)}m_r(q)\right)^{\frac{1}{2r}-\frac{1}{2r^2}}\left(\frac{q}{\phi(q)}\right)^{\frac{1}{r}}.$$
\end{theorem}

\begin{table}[h!]
\begin{center}
  \begin{tabular}{ | c | c | c|c|}
    \hline
    $r$ & $C(r)$ & $D(r)$ & $C^*(r)$ \\ \hline
    2 & 14.275 & 7.844 &  12.364\\ \hline
    3 & 5.135 & 4.390 & 4.849\\ \hline
    4 & 3.556 & 3.289 & 3.448\\ \hline
    5 & 2.878 & 2.740 & 2.823\\ \hline
    6 & 2.496 & 2.410 & 2.463\\ \hline
    7 & 2.248 & 2.189 & 2.227\\ \hline
    8 & 2.075 & 2.030 & 2.060\\ \hline
    9 & 1.945 & 1.911 & 1.934\\ \hline
    $\ge$10 & 1.846 & 1.818 & 1.837\\ \hline
  \end{tabular}
  \caption{Constants in the inequalities for values of $r$ in Theorems \ref{main thm 1} and \ref{main theorem 2}.}\label{Table 1}
\end{center}
\end{table}

The theorem above requires $q$ to be very large, so we have the following theorem which is weaker asymptotically but works for smaller values of $q$.

\begin{theorem}\label{main theorem 2}
    Let $\chi$ be a primitive Dirichlet character modulo $q$. Let $C^{*}(r)$ be defined as in Table \ref{Table 1}. Then, for $q \ge \max\{10^{1143}, 2^{4r-2}\}$, if $r=2$ or $q$ is cubefree, we have
    $$|S_{\chi}(M,N)| \le C^{*}(r)N^{1-\frac{1}{r}} q^{\frac{r+1}{4r^2}}(\log{q})^{\frac{1}{2r}}\left((4r)^{\omega(q)}m_r(q)\right)^{\frac{1}{2r}}\left(\frac{q}{\phi(q)}\right)^{\frac{1}{r}}.$$
\end{theorem}
Theorems~\ref{main thm 1}, \ref{main theorem 2} improve Theorem $A$ in several aspects: it extends the result to $r\geq 3$, allows smaller values of $q$, and removes the hypothesis $N\leq q^{5/8}$. In particular, for $r\leq 949$, the condition $q\geq \max\{10^{1143},2^{4r-2}\}$ reduces to $q\geq 10^{1143}$.

To prove Theorems \ref{main thm 1} and \ref{main theorem 2} we need the following explicit Weil-type inequality, which we consider of independent interest. 

\begin{theorem}[Weil-type inequality]\label{weil}
Let $r\ge 2$ and $q$ be positive integers such that $r=2$ or $q$ is cubefree. Let $\chi$ be a primitive Dirichlet modulo $q$. Let $B \ge 2$ be a real number.  
Then
\begin{equation*}
\sum_{x=1}^q\left|\sum_{1\le b\le B}\chi(x+b)\right|^{2r} \le 2r(4r)^{\omega(q)} B^{2r}m_r(q)\sqrt{q} + \frac{r^{2r}}{r!} B^r q.
\end{equation*}
\end{theorem}

The paper is organized as follows. We prove Theorem \ref{weil} in Section~\ref{sec: weil}. In Section~\ref{General setup} we describe the overall plan for the proofs of our main theorems following the approach for the explicit Burgess inequality detailed in \cite[Theorem 12.6]{IK2004}. We prove other technical lemmas in Section~\ref{sec: technical}, some important bounds in Section~\ref{section asymptotics}, and we prove our main theorems in Section~\ref{sec: main theorems}.

\subsection*{Notation}
Most of the notation in the paper is standard, for example, $\phi(n)$ is Euler's totient function, for a finite set $S$, $\#S$ is the number of elements of $S$, $\lfloor x\rfloor$ is the floor of $x$, for a pair of integers $a,b$, $(a,b) = \gcd(a,b)$. We also use Landau's $O$ notation, i.e., $f(x) = O(g(x))$ if there exists a positive constant $C$ such that $|f(x)| \le C|g(x)|$ for large enough $x$. Some notation we use that is less well-known is $O^*$, where we say $f(x) = O^*(g(x))$ if for all $x\ge x_0$, $|f(x)| \le |g(x)|$. Finally, in some instances to simplify algebraic terms, we use $\phi^* =\frac{\phi(q)}{q}$.

\section{Explicit Weil-type inequality}\label{sec: weil}
This section can be read independently from the rest of the paper. Let $r\ge 2$ be an integer, $B\ge 2$ be a real number, and \begin{align}\label{meet fake mrq}
    m_r(q) = \min \biggl\{ \tau_{2r}(q), \biggl( \frac{\tau(q)}2 \biggr)^{2r-1} , \frac{q}{2r} \biggr\}.
\end{align} 
Furthermore, given integers $b_1, b_2, \dots,b_{2r}$, we define 
$$A_j= \displaystyle\prod_{\substack{1\le i\le 2r\\i\ne j}}\left(b_i - b_j\right).$$

Burgess \cite[Lemma 7]{burgess:1962}, \cite[Lemma 8]{Bur1963} proved the following analogue of Weil's inequality that holds for composite moduli:

\begin{lemma}[Burgess]\label{burgess weil}
 Let $b_1,b_2,\ldots,b_{2r}$ be integers such that at least $r+1$ of them are distinct. Let $f_1(x) = (x-b_1)(x-b_2)\cdots(x-b_r)$ and $f_2(x) = (x-b_{r+1})(x-b_{r+2})\cdots (x-b_{2r})$. Let $\chi$ be a primitive Dirichlet character modulo $q$. If $r=2$ or $q$ is cubefree, then there exists $1\le j\le 2r$ with $A_j\ne 0$ such that
 $$\left|\sum_{x\bmod{q}}\chi(f_1(x)\cdot f_2^{(\phi(q)-1)}(x))\right|\le (4r)^{\omega(q)}\sqrt{q}(A_j,q),$$
 where $(A_j,q)$ is the greatest common divisor of $A_j$ and $q$.
 In particular, 
 $$\left|\sum_{x\bmod{q}}\chi(f_1(x)\cdot f_2^{(\phi(q)-1)}(x))\right|\le \sum_{\substack{j=1\\A_j\ne 0}}^{2r}(4r)^{\omega(q)}\sqrt{q}(A_j,q).$$
 
\end{lemma}

\begin{definition} \label{def sq(r)}
Define
    \begin{align*}
      s_q(r,B) = \sum_{1\le b_1\le B}\sum_{1\le b_2\le B}\cdots\sum_{1\le b_{2r}\le B}
      \min\bigl\{ (A_j, q)\colon 1\le j\le 2r,\, A_j\ne0 \bigr\}.
    \end{align*}
First, observe that from the trivial bound
\begin{align*}
     s_q(r,B)\le \sum_{1\le b_1\le B}\sum_{1\le b_2\le B}\cdots\sum_{1\le b_{2r}\le B}
      \sum_{\substack{1\le j\le 2r\\A_j\ne0}}\left(A_j, q\right)\le B^{2r}q.
\end{align*}
\end{definition}

Now, the following lemmas give complementary upper bounds for this expression.
\begin{lemma}
\label{lem no keep gcd}
Let $s_q(r, B)$ be defined as in Definition \ref{def sq(r)}. Then
$$s_q(r, B)\le 2r \biggl(\frac{\tau(q)}{2}\biggr)^{2r-1}B^{2r}.$$ 
\end{lemma}

\begin{proof}
    Using the bound
    \begin{equation*}
    (A_j,q) = \biggl( \prod_{\substack{1\le i\le 2r \\ i\ne j}} (b_i-b_j), q \biggr) \le \prod_{\substack{1\le i\le 2r \\ i\ne j}} (b_i-b_j,q),
    \end{equation*}
    and exchanging the order of summation, we have
    \begin{equation}
    \label{eq no keep gcd}
    \begin{split}
        s_q(r, B)\le\sum_{1\le j\le 2r}\sum_{\substack{1\le b_1,b_2,\dots,b_{2r}\le B \\ A_j\ne0}} 
      \left(A_j, q\right)
      &\le 
      \sum_{1\le j\le 2r} \sum_{\substack{1\le b_1,b_2,\dots,b_{2r}\le B \\ A_j\ne0}} \prod_{\substack{1\le i\le 2r \\ i\ne j}} (b_i-b_j,q)
      \\
      &= \sum_{1\le j\le 2r} \sum_{1\le b_j\le B} \sum_{\substack{1\le b_1,\dots,b_{j-1},b_{j+1},\dots,b_{2r}\le B \\ b_j\notin\{b_1,\dots,b_{j-1},b_{j+1},\dots,b_{2r}\}}} \prod_{\substack{1\le i\le 2r \\ i\ne j}} (b_i-b_j,q) \\
&= \sum_{1\le j\le 2r} \sum_{1\le b_j\le B} \biggl( \sum_{\substack{1\le b\le B \\ b\ne b_j}} (b-b_j,q) \biggr)^{2r-1} \\
&= 2r \sum_{1\le a\le B} \biggl( \sum_{\substack{1\le b\le B \\ b\ne a}} (b-a,q) \biggr)^{2r-1} .
    \end{split}
    \end{equation}
Using the identity $\sum_{d\mid n} \phi(d) = n$, we write the inner sum as
\[
\sum_{\substack{1\le b\le B \\ b\ne a}} (b-a,q) = \sum_{\substack{1\le b\le B \\ b\ne a}} \sum_{d\mid(b-a,q)} \phi(d) = \sum_{\substack{d\mid q \\ d\le B-1}} \phi(d) \sum_{\substack{1\le b\le B \\ b\ne a \\ d\mid b-a}} 1 \le \sum_{\substack{d\mid q \\ d\le B}} \phi(d) \frac {\lfloor B\rfloor}d = \lfloor B\rfloor \sum_{\substack{d\mid q \\ d\le B}} \frac{\phi(d)}d,
\]
where since the divisors of $q$ (other than $\sqrt{q}$ when $q$ is a square) pair up on either side of $\sqrt q$ and $B<\sqrt q$,
\begin{equation}
\label{eq tau(q) upper bound}
    \lfloor B\rfloor\sum_{\substack{d\mid q \\ d\le B}} \frac{\phi(d)}d 
    \le 
    \lfloor B\rfloor \sum_{\substack{d\mid q \\ d\le B}} 1 \le \frac{\lfloor B\rfloor}2\tau(q).
\end{equation}
Putting inequalities \eqref{eq no keep gcd} and \eqref{eq tau(q) upper bound} together,
\begin{align*} 
s_q(r, B) &\le 2r \sum_{1\le a\le B} \biggl( \sum_{\substack{1\le b\le B \\ b\ne a}} (b-a,q) \biggr)^{2r-1} \\
&\le 2r \sum_{1\le a\le B} \biggl(\frac{\lfloor B\rfloor}2\tau(q) \biggr)^{2r-1} = 2r \left(\frac{\tau(q)}{2}\right)^{2r-1}\lfloor B\rfloor^{2r}, 
\end{align*}
as stated in the lemma.
\end{proof}

We bound the sum $s_q(r, B)$ below without using the initial inequality on the $\gcd$ in Lemma~\ref{lem no keep gcd}.

\begin{lemma}\label{lem keep gcd}
Let $s_q(r, B)$ be defined as in Definition \ref{def sq(r)}. Then
$$s_q(r, B)\le 2r B^{2r} \tau_{2r}(q).$$
\end{lemma}

\begin{proof}
Exchanging the sum as before and using the identity $\sum_{d\mid n} \phi(d) = n$,
\begin{align*}
s_q(r, B) \le \sum_{1\le j\le 2r} \mathop{\sum\nolimits}\limits_{\substack{1\le b_1,b_2,\dots,b_{2r}\le B \\ A_j\ne0}} (A_j,q) &= \sum_{1\le j\le 2r} \mathop{\sum\nolimits}\limits_{\substack{1\le b_1,b_2,\dots,b_{2r}\le B \\ A_j\ne0}} \sum_{d\mid (A_j,q)} \phi(d) \\
&\le \sum_{1\le j\le 2r} \sum_{\substack{d\mid q \\ d\le B^{2r-1}}} \phi(d) \sum_{\substack{1\le b_1,b_2,\dots,b_{2r}\le B \\ A_j\ne0 \\ d\mid A_j}} 1.
\end{align*}
Let $d_1\cdots d_{2r}$ be an ordered factorization of $d$ with $d_j=1$, then $d\mid A_j$ if $d_i\mid(b_i-b_j)$ for each $1\le i\le 2r$. Every inner summand arises in this way for at least one such ordered factorization, therefore
\begin{align*}
\sum_{1\le j\le 2r} \sum_{\substack{d\mid q \\ d\le B^{2r-1}}} \phi(d) \sum_{\substack{1\le b_1,b_2,\dots,b_{2r}\le B \\ A_j\ne0 \\ d\mid A_j}} 1 \le \sum_{1\le j\le 2r} \sum_{\substack{d\mid q \\ d\le B^{2r-1}}} \phi(d) \sum_{1\le b_j\le B} \sum_{\substack{d_1\cdots d_{2r}=d \\ d_j=1}} \sum_{\substack{1\le b_1\le B \\ b_1\ne b_j \\ d_1\mid(b_1-b_j)}} \cdots \sum_{\substack{1\le b_{2r}\le B \\ b_{2r}\ne b_j \\ d_{2r}\mid(b_{2r}-b_j)}} 1.
\end{align*}
Note that there are $2r-1$ inner summands with~$b_j$ fixed. Bounding each inner summand by $\displaystyle \frac{\lfloor B\rfloor}{d_i}$,
\begin{align}\notag
\label{eq keep gcd}
s_q(r, B)&\le \sum_{1\le j\le 2r} \sum_{\substack{d\mid q \\ d\le B^{2r-1}}} \phi(d) \sum_{1\le b_j\le B}  \sum_{\substack{d_1\cdots d_{2r}=d \\ d_j=1}} \sum_{\substack{1\le b_1\le B \\ b_1\ne b_j \\ d_1\mid(b_1-b_j)}} \cdots \sum_{\substack{1\le b_{2r}\le B \\ b_{2r}\ne b_j \\ d_{2r}\mid(b_{2r}-b_j)}} 1 \\
&\le \sum_{1\le j\le 2r} \sum_{\substack{d\mid q \\ d\le B^{2r-1}}} \phi(d) \sum_{1\le b_j\le B} \sum_{\substack{d_1\cdots d_{2r}=d \\ d_j=1}} \prod_{\substack{1\le i\le 2r \\ i\ne j}} \frac {\lfloor B\rfloor}{d_i} \\ \notag
&= \lfloor B\rfloor^{2r} \sum_{1\le j\le 2r} \sum_{\substack{d\mid q \\ d\le B^{2r-1}}} \frac{\phi(d)}d \sum_{\substack{d_1\cdots d_{2r}=d \\ d_j=1}} 1 
= 2r \lfloor B\rfloor^{2r} \sum_{\substack{d\mid q \\ d\le B^{2r-1}}} \frac{\phi(d)}d \tau_{2r-1}(d).
\end{align}
Since
\[
\sum_{d\mid q} \frac{\phi(d)}d \tau_{2r-1}(d) = 
\prod_{p^k\|q} \biggl( 1 + \frac{p-1}p \bigl( \tau_{2r}(p^k) - 1 \bigr) \biggr) \le \tau_{2r}(q),
\]
combining with inequality \eqref{eq keep gcd}, we obtain
\[
s_q(r, B)
\le
\sum_{1\le j\le 2r} \mathop{\sum\nolimits}\limits_{\substack{1\le b_1,b_2,\dots,b_{2r}\le B \\ A_j\ne0}} (A_j,q) \le 2r \lfloor B\rfloor^{2r} \tau_{2r}(q),
\]
as stated above.
\end{proof}

Using the lemmas above, we prove the explicit Weil-type inequality.

\begin{proof}[Proof of Theorem \ref{weil}]
As done in the proof of Theorem 1.1 in \cite{ET1}, letting $f_1(x)$ and $f_2(x)$ be defined as in Lemma \ref{burgess weil}, we have
\begin{equation*}\label{poly}
\sum_{x=1}^q\left|\sum_{b\in\mathcal{B}}\chi(x+b)\right|^{2r}= \sum_{1\le b_1,b_2,\ldots,b_{2r}\le B}\sum_{x\bmod{q}}\chi(f_1(x)f_2^{(\phi(q)-1)}(x)).
\end{equation*}
Let $(b_1,b_2,\ldots,b_{2r})$ be a good tuple if at least $r+1$ of them are distinct. Let $\cB_1$ be the set of good tuples and $\cB_2$ be the rest. By Lemma \ref{burgess weil},
$$\sum_{(b_1,\ldots,b_{2r})\in \cB_1}\sum_{x\bmod{q}}\chi(f_1(x)f_2^{(\phi(q)-1)}) \le (4r)^{\omega(q)}\sqrt{q} s_q(r, B).$$
Now using Lemmas \ref{lem no keep gcd} and \ref{lem keep gcd}, we obtain the upper bound
\begin{align}\label{eq good tuples}
  \sum_{(b_1,\ldots,b_{2r})\in \cB_1}\sum_{x\bmod{q}}\chi(f_1(x)f_2^{(\phi(q)-1)})\le 2r(4r)^{\omega(q)}m_r(q)\sqrt{q}\lfloor B\rfloor^{2r},  
\end{align}
where $m_r(q)$ is given in \eqref{meet fake mrq}.

For the other tuples $\cB_2$, we use the trivial bound,  
$$\sum_{x\bmod{q}} \chi(f_1(x)f_2^{(\phi(q)-1)}(x)) \le q.$$
If $(b_1, b_2,\ldots, b_{2r}) \in \cB_2$, then it has at most $r$ distinct values. There are $\binom{\lfloor B\rfloor}{r}$ ways of choosing the $r$ values and each $b_i$ has $r$ choices to make, therefore, 
\begin{align}\label{eq bad tuples}
    \sum_{(b_1,b_2,\ldots,b_{2r})\in \cB_2}\sum_{x\bmod{q}}\chi(f_1(x)f_2^{(\phi(q)-1)}(x)) \le r^{2r} \binom{\lfloor B\rfloor}{r} q \le \frac{r^{2r} \lfloor B\rfloor^r}{r!}q.
\end{align}
Combining the bounds in \eqref{eq good tuples} and \eqref{eq bad tuples}, we get the desired upper bound in the Weil-type inequality.
\end{proof}

\section{General setup}\label{General setup}

Let $r\ge 2, N$ be positive integers, and let $q$ be a positive integer that is cubefree when $r\ge 3$ for the remainder of the paper.
Given a primitive Dirichlet character $\chi$ modulo $q$ and integers~$M$ and ~$N\ge1$, we consider the character sum
\begin{equation*}
    S_\chi\left(M, N\right)=\sum_{M<n \le M+N}\chi(n).
\end{equation*}
Given an integer $r\ge2$, we use induction to prove a bound of the form
\begin{equation} \label{eq desired bound}
| S_\chi(M,N) | \le E_{q,r}(N), \quad\text{where }
E_{q,r}(N) = C N^{1-1/r} q^{\frac{r+1}{4r^2}} (\log q)^{\frac{1}{2r}} T(q),
\end{equation}
where $C\ge 1$ is a constant to be determined later, $m_r(q)$ is defined in \eqref{meet Mr. Q}, and
\begin{equation}\label{meet T(q)}
T(q) = \begin{cases}\bigl((4r)^{\omega(q)} m_r(q)\bigr)^{\left({\frac{1}{2r}-\frac{1}{2r^2}}\right)} \bigl( \frac q{\phi(q)} \bigr)^{1/r}& \text{ under the conditions of Theorem } \ref{main thm 1},\\
\bigl((4r)^{\omega(q)} m_r(q)\bigr)^{\left({\frac{1}{2r}}\right)} \bigl( \frac q{\phi(q)} \bigr)^{1/r}& \text{ under the conditions of Theorem } \ref{main theorem 2} 
.\end{cases}
\end{equation}

We take $1\le N\le q^{\frac{1}{4}+\frac{1}{4r}}$ as the base case. In this case, the inequality~\eqref{eq desired bound} follows from the trivial bound $|S_\chi(M,N)| \le N$, since
\begin{align*}
&N \le C N^{1-1/r} q^{\frac{r+1}{4r^2}} (\log q)^{\frac{1}{2r}} T(q) \\
\implies &N \le q^{\frac{r+1}{4r}} \bigl( C (\log q)^{\frac{1}{2r}} T(q) \bigr)^r,
\end{align*}
when $N \le q^{\frac{1}{4}+\frac{1}{4r}}$ (because $C(\log q)^{\frac{1}{2r}} T(q) \ge 1$).
\\
Therefore, in the induction step, we assume
\begin{align}\label{inductive hypothesis}
\text{for $N > q^{\frac{1}{4}+\frac{1}{4r}}$, the bound in \eqref{eq desired bound} holds for all smaller positive integers.}
\end{align}

We use the following notation for the remainder of this paper.
\begin{definition} \label{def A B v}
For some real numbers $A,B\ge2$, let
\begin{itemize}
\item $\cA=\{1\le a \le A\colon a\in\mathbb{Z},\, (a, q)=1\}$,
\item $\cB=\{1\le b\le B\colon b\in\mathbb{Z}\}$,
\item $v_\cA(x)=\#\{a\in\cA,\, n\in(M,M+N] \colon \bar{a}n\equiv x \pmod{q}\}$, where $\bar{a}$ denotes the multiplicative inverse of~$a$ modulo~$q$.
\end{itemize}
\end{definition}

We prove the following lemma using the notations in \eqref{eq desired bound}, \eqref{meet T(q)}, and Definition \ref{def A B v}.
\begin{lemma}\label{eq after Holder}
For positive integers $r\ge 2$ and $q$, where $q$ is cubefree for $r\ge 3$, we have, under the inductive hypothesis in \eqref{inductive hypothesis}, 
    \begin{align*}
    \left|  S_\chi\left(M, N\right)\right|
    &\le
    \frac{1}{\#\cA\cdot \#\cB}\left(\sum_{x=1}^q v_\cA(x)\right)^{1-\frac{1}{r}}
    \left(\sum_{x=1}^q v^2_{\cA}(x)\right)^{\frac{1}{2r}}\left(\sum_{x=1}^q\left|\sum_{b \in \cB}\chi\left(x+b\right)\right|^{2r}\right)^\frac{1}{2r}
    \\
    &\qquad{}+
\frac{2C}{2-1/r} A^{1-1/r} \frac{(\lfloor B\rfloor+1)^{2-1/r}}{\#\cB} q^{\frac{r+1}{4r^2}} (\log q)^{\frac{1}{2r}} T(q),   
\end{align*}
for some constant $C\ge 1$.
\end{lemma}

\begin{proof}
Under the inductive hypothesis in \eqref{inductive hypothesis}, we have,
by equation~(28) in \cite{ET2},
\begin{align}\label{eq 28 in ET16}
      \left|  S_\chi\left(M, N\right)\right| 
       \le 
       \frac{1}{\#\cA\cdot \#\cB}\sum_{x=1}^{q}v_\cA(x)\left|\sum_{1\le b \le B}\chi(x+b)\right| +
    \frac{1}{\#\cA\cdot \#\cB}\sum_{a\in \cA}\sum_{b\in \cB}2E_{q,r}(ab),
    \end{align}
where $E_{q,r}(ab)$ is defined as in \eqref{eq desired bound}.

Since $E_{q,r}(N)$ is an increasing function of~$N$,
\[
\frac{1}{\#\cA\cdot \#\cB}\sum_{a\in \cA}\sum_{b\in \cB}2E_{q,r}(ab) \le \frac{1}{\#\cA\cdot \#\cB}\sum_{a\in \cA}\sum_{b\in \cB}2E_{q,r}(Ab) = \frac{2}{\#\cB}\sum_{b\in \cB} E_{q,r}(Ab),
\]
and by definition,
\begin{align*}
\frac{2}{\#\cB}\sum_{b\in \cB} E_{q,r}(Ab) &= \frac{2}{\#\cB}\sum_{b\in \cB} C (Ab)^{1-1/r} q^{\frac{r+1}{4r^2}} (\log q)^{\frac{1}{2r}} T(q) \\
&= \frac{2C}{\#\cB} A^{1-1/r} q^{\frac{r+1}{4r^2}} (\log q)^{\frac{1}{2r}} T(q) \sum_{b\in \cB} b^{1-1/r}.
\end{align*}
Since
\[
\sum_{b\in \cB} b^{1-1/r} \le \int_1^{\lfloor B\rfloor+1} t^{1-1/r} \,{\rm d}t = \frac{(\lfloor B\rfloor+1)^{2-1/r}-1}{2-1/r} < \frac{(\lfloor B\rfloor+1)^{2-1/r}}{2-1/r},
\]
we conclude that
\begin{equation} \label{eq E sum bound}
\frac{1}{\#\cA\cdot \#\cB}\sum_{a\in \cA}\sum_{b\in \cB}2E_{q,r}(ab) < \frac{2C}{2-1/r} A^{1-1/r} \frac{(\lfloor B\rfloor+1)^{2-1/r}}{\#\cB} q^{\frac{r+1}{4r^2}} (\log q)^{\frac{1}{2r}} T(q).
\end{equation}

Therefore, combining \eqref{eq 28 in ET16}, \eqref{eq E sum bound}, and applying H{\"o}lder's inequality, we get
\begin{align*}
    \left|  S_\chi\left(M, N\right)\right|
    &\le
    \frac{1}{\#\cA\cdot \#\cB}\left(\sum_{x=1}^q v_\cA(x)\right)^{1-\frac{1}{r}}
    \left(\sum_{x=1}^q v^2_{\cA}(x)\right)^{\frac{1}{2r}}\left(\sum_{x=1}^q\left|\sum_{b \in \cB}\chi\left(x+b\right)\right|^{2r}\right)^\frac{1}{2r}
    \\
    &\qquad{}+
\frac{2C}{2-1/r} A^{1-1/r} \frac{(\lfloor B\rfloor+1)^{2-1/r}}{\#\cB} q^{\frac{r+1}{4r^2}} (\log q)^{\frac{1}{2r}} T(q),   
\end{align*}
as stated in the lemma.
\end{proof}
In the next section, we handle the two sums involving $v_{\cA}$.

\section{The $v_{\cA}$ sums}\label{sec: technical}

Recall from Definition~\ref{def A B v} that for any real number $A\ge 2$, $\cA=\{1\le a \le A\colon a\in\mathbb{Z},\, (a, q)=1\}$ and $v_\cA(x)=\#\{a\in\cA,\, n\in(M,M+N] \colon \bar{a}n\equiv x \pmod{q}\}$, where $\bar{a}$ denotes the multiplicative inverse of~$a$ modulo~$q$. We first prove the following lemma.

\begin{lemma} \label{v_A sum}
Let $\cA$ and $v_\cA(x)$ be defined as in Definition \ref{def A B v}. Then 
$$\displaystyle \sum_{x\pmod q} v_\cA(x) = \#\cA\cdot N.$$
\end{lemma}

\begin{proof}
By definition,
\[
\sum_{x\pmod q} v_\cA(x) = \sum_{x\pmod q} \#\bigl\{ a\in\cA,\, n\in(M,M+N] \colon \bar{a}n\equiv x \pmod{q} \bigr\}.
\]
Summing over all $x$ satisfying $x\equiv\bar an\pmod q$, each ordered pair $(a,n)$ with $a\in\cA$ and $n\in(M,M+N]$ is counted exactly once. Therefore, we have the equality stated in the lemma.
\end{proof}

\begin{lemma} \label{JKL page 8}
Let $\cA$ and $v_\cA(x)$ be defined as in Definition \ref{def A B v}. Then
$$\displaystyle \sum_{x\pmod q} v_\cA(x)^2 \le \sum_{a_1\in\cA}
\sum_{a_2\in\cA} \biggl( 1 + \frac{N(a_1,a_2)}{\max\{a_1,a_2\}} \biggr).$$
\end{lemma}

\begin{proof}
This is the first displayed inequality under Equation (4.5)
in \cite{JKL}.
\end{proof}

The following Lemma from \cite[Lemma 16]{BrokerSutherland2010} allows us to bound the sum of $\phi(n)/n^2$.\begin{lemma}\label{phi sum lemma}
For all $N \geq 1$, we have
\[
\sum_{n\leq N} \frac{\phi(n)}{n^2}
\leq
\frac{6}{\pi^2}\log N
+\delta
+\frac{\log N+1}{3N},
\]
where
\[
\delta
=
\frac{6\gamma}{\pi^2}
-\frac{36\zeta'(2)}{\pi^4}
\approx 0.6974.
\]
\end{lemma}


The following proposition asymptotically improves the inequality used in \cite{JKL} by a factor of $6/\pi^2$.

\begin{proposition}
\label{real v_A square sum}
Let $A\ge 2$ be a real number, $N\ge2$ be an integer, and $\cA, v_\cA(x)$ be defined as in Definition \ref{def A B v}. Then
$$\displaystyle \sum_{x\pmod q} v_\cA(x)^2 \le \#\cA^2 + 2AN \biggl( \frac6{\pi^2} \log A + \delta + \frac{\log A+1}{3A} \biggr),$$ 
where $\delta$ is defined in Lemma~\ref{phi sum lemma}.
\end{proposition}

\begin{proof}
Starting from Lemma~\ref{JKL page 8},
\begin{align*}
\sum_{x\pmod q} v_\cA(x)^2 \le \sum_{a_1\in\cA} \sum_{a_2\in\cA} \biggl( 1 + \frac{N(a_1,a_2)}{\max\{a_1,a_2\}} \biggr) &= (\#\cA)^2 + N \sum_{a_1\in\cA} \sum_{a_2\in\cA} \frac{(a_1,a_2)}{\max\{a_1,a_2\}} \\
&\le (\#\cA)^2 + N \sum_{1\le a_1\le A} \sum_{1\le a_2\le A} \frac{(a_1,a_2)}{\max\{a_1,a_2\}}.
\end{align*}
In the remaining sum we use our beloved trick $g = \sum_{d\mid g} \phi(d)$ to obtain
\begin{align*}
\sum_{1\le a_1\le A} \sum_{1\le a_2\le A} \frac{(a_1,a_2)}{\max\{a_1,a_2\}} &= \sum_{a_1\in\cA} \sum_{a_2\in\cA} \frac1{\max\{a_1,a_2\}} \sum_{d\mid(a_1,a_2)} \phi(d) \\
&= \sum_{d\le A} \phi(d) \sum_{\substack{1\le a_1\le A \\ d\mid a_1}} \sum_{\substack{1\le a_2\le A \\ d\mid a_2}} \frac1{\max\{a_1,a_2\}} \\
&= \sum_{d\le A} \frac{\phi(d)}d \sum_{1\le b_1\le A/d} \sum_{1\le b_2\le A/d} \frac1{\max\{b_1,b_2\}} \\
&\le 2 \sum_{d\le A} \frac{\phi(d)}d \sum_{1\le b_1\le A/d} \sum_{1\le b_2\le b_1} \frac1{b_1} \\
&= 2 \sum_{d\le A} \frac{\phi(d)}d \sum_{1\le b_1\le A/d} 1 \le 2A \sum_{d\le A} \frac{\phi(d)}{d^2}.
\end{align*}
Now the proposition follows from Lemma~\ref{phi sum lemma}.
\end{proof}




\section{Important bounds}\label{section asymptotics}
As we can see in Lemma \ref{eq after Holder}, we need to make some choices for $A$ and $B$ to get a Burgess inequality. The following lemmas guarantee that the $A$ and $B$ we pick later satisfy certain properties that are required in the proofs of Theorem \ref{main thm 1} and Theorem \ref{main theorem 2}.

\begin{lemma} \label{B will be at least 10 if it kills us lemma}
For $r\ge 2$, let $a(r) = 2\left(3.0758r+1.38402\log(4r)-1.5379\right)\log{2}$. If $q \ge e^{e^{a(r)}}$, then $$B_1(r,q)\ge 2^{2-2/r}r^2\left(\frac{1-\frac{1}{r}}{r!}\right)^{\frac{1}{r}},$$
where 
\begin{equation}\label{def B_1}
B_1(r,q) = r^2q^{\frac{1}{2r}}\biggl( \frac{r-1}{r!2r(4r)^{\omega(q)}m_r(q)} \biggr)^{1/r}.
\end{equation}
\end{lemma}
 
\begin{proof} 
Since $m_r(q) \le \bigl( \tau(q)/2 \bigr)^{2r-1}$ by definition~\eqref{meet Mr. Q} 
we have
\[
B_1(r,q) \ge r^2q^{\frac{1}{2r}}\biggl( \frac{2^{2r-1}(r-1)}{2r (r!)(4r)^{\omega(q)}(\tau(q))^{2r-1}} \biggr)^{1/r}. 
\]
We have $2^{\omega(q)} \le 2^{1.38402\log q/\log\log q}$ by~\cite{Nicolas-Robin1983} and $\tau(q) \le 2^{1.5379\log q/\log\log q}$ \cite[Theorem 1]{Nicolas-Robin1983}, therefore
$$(4r)^{\omega(q)}\le q^{\frac{1.38402\log{2} \log{(4r)}}{\log\log{q}}},$$
and
\[
B_1(r,q) \ge 2^{2-2/r}r^2\left(\frac{1-\frac{1}{r}}{r!}\right)^{\frac{1}{r}} q^{\frac{1}{2r}-\frac{1.38402(\log{2})\log{(4r)}}{r \log\log{q}}-\frac{1.5379 (2r-1) \log{2}}{r \log\log{q}}}. 
\]

When $q\ge e^{e^{a(r)}}$, we have
$$\frac{1}{2r} - \frac{1.38402(\log{2})\log{(4r)}}{r \log{\log{q}}}-\frac{1.5379(2r-1)\log{2}}{r \log{\log{q}}} \ge 0,$$
so
$$B_1(r,q) \ge 2^{2-2/r}r^2\left(\frac{1-\frac{1}{r}}{r!}\right)^{\frac{1}{r}}$$
as desired.
\end{proof}
\begin{remark}
The doubly-exponential lower bound on \(q\) comes from the term
\(\log\log q\) appearing in the estimate
$$
2^{\omega(q)}
\le
2^{1.38402\,\frac{\log q}{\log\log q}}.
$$
Indeed, by the definition of \(a(r)\), the condition needed in the proof is
ensured by $
\log\log q \ge a(r),$ which leads to $q \ge e^{e^{a(r)}}.$
\end{remark}

\begin{lemma} \label{B will be at least 10 if it kills us lemma 2}
Let $r\ge 2$. If $q \ge 2^{4r-2}$, then $$B_2(r,q)\ge 2^{2-2/r}r^2\left(\frac{1-\frac{1}{r}}{r!}\right)^{\frac{1}{r}},$$
where 
\begin{equation}\label{def B_2}B_2(r,q) = r^2q^{\frac{1}{2r}}\biggl( \frac{r-1}{r!2r} \biggr)^{1/r}.
\end{equation}
\end{lemma}
\begin{proof}
    Indeed, for $q \ge 2^{4r-2}$,
    $$r^2q^{\frac{1}{2r}}\biggl( \frac{r-1}{r!2r} \biggr)^{1/r} \ge 2^{2-2/r}r^2\left(\frac{1-\frac{1}{r}}{r!}\right)^{\frac{1}{r}}.$$
\end{proof}

\begin{lemma} \label{B upper bound lemma}
For integers $r, q\ge2$, we have 
$$B_i(r,q) \le erq^{\frac{1}{2r}},$$ 
for $i = 1,2$ with the $B_i(r, q)$ defined as in \eqref{def B_1} and \eqref{def B_2}. 
\end{lemma}

\begin{proof}
Using $\frac{r-1}{2r (4r)^{\omega(q)} m_r(q)} \le 1$ and $r! \ge \left(\frac{r}{e}\right)^r$,\footnote{Since $e^n = \sum_{k=0}^\infty \frac{n^k}{k!} \ge \sum_{k=n}^n \frac{n^k}{k!} = \frac{n^n}{n!}$, we have $n! \ge \frac{n^n}{e^n}$.} we get
$$\frac{B_1(r,q)}{r^2q^{\frac{1}{2r}}} = \left(\frac{r-1}{r! (2r)(4r)^{\omega(q)}m_r(q)}\right)^{\frac{1}{r}} \le \frac{e}{r}.$$
Using $r! \ge \left(\frac{r}{e}\right)^r$, we get
$$\frac{B_2(q,r)}{q^{\frac{1}{2r}}r^2} = \left(\frac{r-1}{r!(2r)}\right)^{\frac{1}{r}} \le \frac{e}{r}.$$

\end{proof}

To prove our main theorems, we also need some lower bounds on $A$. In particular, we need $A \ge 31$, and $A\ge 2^{\omega(q)} q /\phi(q)$. The following lemmas cover this for two different choices of $A$.

\begin{lemma} \label{A lower bound lemma}
Let $N \ge q^{\frac{1}{4}+\frac{1}{4r}}$ and $r\ge 2$ be positive integers. Suppose 
$$A_i(r,q) = \frac{\kappa N}{B_i(r,q)},$$
for $i=1,2$ where $B_1(r,q)$ is defined as in \eqref{def B_1}, $B_2(r,q)$ as in \eqref{def B_2}, and let
\begin{equation} \label{choice of kappa}
\kappa = \biggl( \frac{(B_i(r,q)-1) B_i(r,q)^{1-1/r} }{2 r (B_i(r,q)+1)^{2-1/r}}
\biggr)^{r/(r-1)}.
\end{equation}
Let $t_1(r) = e^{e^{a(r)}}$, where $a(r)$ is defined as in Lemma \ref{B will be at least 10 if it kills us lemma} and $t_2(r) = 2^{4r-2}$. Then, for $q\ge \max\bigl\{ 10^{1143}, (10r)^{18}, t_i(r) \bigr\}$, we have 
$$A_i(r,q) \ge \max\biggl\{31, \frac{2^{\omega(q)}q}{\phi(q)}\biggl\}.$$
\end{lemma}

\begin{proof}
Since the proof is the same for $i\in \{1, 2\}$, let $B = B_i(r,q)$ and $A = A_i(r,q)$ for brevity. By Lemmas  \ref{B will be at least 10 if it kills us lemma} and \ref{B will be at least 10 if it kills us lemma 2}, we have $B\ge 4$ for $q\ge t_i(r)$.   
Therefore, for $q\ge t_i(r)$,
$$\kappa \ge \frac{36}{125}\left(\frac{1}{2r}\right)^{r/(r-1)}.$$

Then, by Lemma~\ref{B upper bound lemma},
\[
A = \frac{\kappa N}B \ge \frac{\kappa q^{\frac{1}{4}+\frac{1}{4r}}}{erq^{\frac{1}{2r}}} \ge \frac{36}{125e}\frac{q^{\frac{1}{4}-\frac{1}{4r}}}{r(2r)^{r/(r-1)}}.
\]
The right-hand side exceeds~$31$ once
\[
q \ge \left( 31\left(\frac{125e}{36}\right)r(2r)^{r/(r-1)} \right)^{4r/(r-1)},
\]
and one can check that this right-hand side is at most $\max\{10^{32},(25r)^8\}$ for all real $r\ge2$.

Moreover, we have $2^{\omega(q)} \le 2^{1.38402\log q/\log\log q}$ by~\cite{Nicolas-Robin1983} and for $q\ge 10^{13}$
\[
\frac q{\phi(q)} < e^\gamma \log\log q + \frac{2.50637}{\log\log q} = \left(e^{\gamma}+\frac{2.50637}{(\log{\log{q}})^2}\right)\log\log q < 2\log\log{q}
\]
by~\cite[Theorem~15]{RS1962}. Thus it suffices to show that
\[
\frac{q^{\frac{1}{4}-\frac{1}{4r}}}{r(2r)^{r/(r-1)}} \ge 2^{1.38402\log q/\log\log q}\cdot \left(e^{\gamma}+\frac{2.50637}{(\log{\log{q}})^2}\right)\log\log q \cdot \frac{125e}{36}
\]
or equivalently
\[
q \ge \left( r(2r)^{r/(r-1)} \cdot 2^{1.38402\log q/\log\log q}\cdot \left(e^{\gamma}+\frac{2.50637}{(\log{\log{q}})^2}\right)\log\log q \cdot\left(\frac{125e}{36}\right)\right)^{(\frac{1}{4}-\frac{1}{4r})^{-1}}.
\]
When $2\le r\le 9$ we check computationally that this inequality holds for $q\ge10^{1143}$. For $r\ge10$, we use the checkable bound
\[
\bigl( r(2r)^{r/(r-1)} \bigr)^{(\frac{1}{4}-\frac{1}{4r})^{-1}} \le 741r^8
\]
together with the trivial $(\frac{1}{4}-\frac{1}{4r})^{-1} \le \frac{40}9$ to see that it suffices for
\[
q \ge 741r^8 \bigl( 2^{1.38402\log q/\log\log q}\cdot 2\log\log q \big)^{40/9}\left(\frac{125e}{36}\right)^{40/9}.
\]
We can computationally show the inequality
\[
2^{1.38402\log q/\log\log q}\cdot 2\log\log q \le q^{1/8}
\]
for $q\ge 10^{1008}$, and so it suffices to have
\[
q \ge 741r^8 \bigl( q^{1/8} \big)^{40/9}\left(\frac{125e}{36}\right)^{40/9},
\]
which holds when $q\ge (10r)^{18} \ge (741r^8)^{9/4}\left(\frac{125e}{36}\right)^{10}$.
\end{proof}

\begin{lemma}\label{lem: alpha bound}
Let $r,q\ge 2$ and $N \ge q^{\frac{1}{4}+\frac{1}{4r}} $ be integers and let $i \in \{1,2\}$. Let $A = A_i(r,q)$, $B = B_i(r,q)$, $\kappa$ be defined as in \eqref{choice of kappa}, and $t_i(r)$ be defined as in Lemma \ref{A lower bound lemma}.
Let
\begin{equation*}
    \alpha
    =
    \left[ \frac{\left(\kappa\left(\frac{\phi(q)}{q}\right)+\frac{B2^{\omega(q)-1}}{N} \right)^2}{B}
    +
    2\kappa\left(\frac{6}{\pi^2}\log{(A)}+\delta+\frac{\log{(A)}+1}{3A}\right)\right],
\end{equation*}
where $\delta = \frac{6\gamma}{\pi^2} - \frac{36\zeta'(2)}{\pi^4} \approx 0.6974 $.
Then, for $q\ge \max\{10^{1143}, (10r)^{18}, t_i(q)\}$,
$$\alpha \le \frac{32}{31}\kappa \log{q}.$$
\end{lemma}

\begin{proof}
Since $q\ge \max\{10^{1143}, (10r)^{18}, t_i(q)\}$, we have $A \ge \frac{2^{\omega(q)}q}{\phi(q)}$ by Lemma \ref{A lower bound lemma}. Using $AB = \kappa N$, we obtain
\begin{align*}
\frac{\left(\kappa\phi^*+\frac{B2^{\omega(q)-1}}{N} \right)^2}{B} &= \frac{A\left(\kappa\phi^*+\frac{\kappa 2^{\omega(q)-1}}{A}\right)^2}{\kappa N} = \frac{\kappa A \left(\phi^* + \frac{2^{\omega(q)-1}}{A}\right)^2}{N}\\
&= \frac{\kappa(A\phi^* + 2^{\omega(q)-1})^2}{AN} \le \frac{\kappa \left(\frac{3}{2} A\phi^*\right)^2}{AN} \le \frac{9\kappa A}{4 N} \le \frac{9\kappa}{4}.
\end{align*}

We also have $A\ge 31$, which implies
$$\delta + \frac{\log{(A)}+1}{3A} \le 0.217\log{A}.$$

We may assume $N\le q^\frac{5}{8}$ because otherwise the P\'olya--Vinogradov inequality produces a better bound (for an explicit P\'olya-Vinogravov inequality see \cite{FroSo}). Thus $A\le q^{\frac{5}{8}}$, and for $q \ge 10^{900}$,
\begin{align*}
\alpha &\le \frac{9\kappa}{4} + 2\kappa\log{A}\left(\frac{6}{\pi^2} + 0.217\right) \le \frac{9\kappa}{4} + \frac{5}{4}\kappa\log{q} \left(\frac{6}{\pi^2}+0.217\right) \\&= \frac{5}{4}\kappa\log{q}\left(\frac{9}{5\log{q}} + \frac{6}{\pi^2} + 0.217\right) < \frac{32}{31}\kappa\log{q}.
\qedhere
\end{align*}
\end{proof}

\section{Proof of Main Theorems}\label{sec: main theorems}
In this section, we prove Theorems \ref{main thm 1} and \ref{main theorem 2}. Recall that we let $r\ge 2$, $q$ be positive integers and $i\in \{1, 2\}$. 
We let $i = 1$ be the case for Theorem \ref{main thm 1}, and $i = 2$ for Theorem \ref{main theorem 2}. 
Let 
\begin{align}\label{eq T1}
    T_1(q) = \bigl((4r)^{\omega(q)} m_r(q)\bigr)^{\left({\frac{1}{2r}-\frac{1}{2r^2}}\right)} \biggl( \frac q{\phi(q)} \biggr)^{1/r},
\end{align}
and
\begin{align}\label{eq T2}
    T_2(q) = \bigl((4r)^{\omega(q)} m_r(q)\bigr)^{\left({\frac{1}{2r}}\right)} \biggl( \frac q{\phi(q)} \biggr)^{1/r}.
\end{align}

Let $t_i(r)$ be defined as in Lemma \ref{A lower bound lemma}. Our goal is to prove that for $q \ge \max\{10^{1143}, t_i(r)\}$, there exists a constant $C(r)$ depending on $r$ (but not $q$) such that
$$|S_{\chi}(M,N)| \le C(r)N^{1-\frac{1}{r}}q^{\frac{r+1}{4r^2}}(\log{q})^{\frac{1}{2r}}T_i(q).$$

\begin{lemma}\label{lemma to set up the main theorems}
For $i\in\{1, 2\}$ and $T_i$ be defined as in \eqref{eq T1} and \eqref{eq T2}. Let $s, \alpha, u, $ and $w$ be defined as in \eqref{eq s}, \eqref{eq alpha}, \eqref{eq u}, and \eqref{w}, respectively. Then
    \begin{equation*}
\begin{split}
    \left| S_\chi(M, N) \right|
    \le
    \left(\frac{B}{B-1}\right)N^{1-1/r}q^{(r+1)/{4r^2}}(\log q)^{\frac{1}{2r}}
    \left(
    u(\alpha w)^{\frac{1}{2r}}
    +sT_i(q)
    \right).
\end{split} 
\end{equation*}
\end{lemma}
\begin{proof}
    Fix $i\in \{1,2\}$ and let $q\ge \max\{10^{1143},t_i(r)\}$. From Lemma \ref{eq after Holder}, once we choose $A$ and $B$ as described in Definition \ref{def A B v}, replacing $T(q)$ with $T_i(q)$ and $C(r)$ with $C$, we have
\begin{align*} 
    \left|  S_\chi\left(M, N\right)\right|
    &\le
    \frac{1}{\#\cA\cdot \#\cB}\left(\sum_{x=1}^q v_\cA(x)\right)^{1-\frac{1}{r}}
    \left(\sum_{x=1}^q v^2_{\cA}(x)\right)^{\frac{1}{2r}}\left(\sum_{x=1}^q\left|\sum_{b \in \cB}\chi\left(x+b\right)\right|^{2r}\right)^\frac{1}{2r}
    \\
    &\qquad{}+
\frac{2C}{2-1/r} A^{1-1/r} \frac{(\lfloor B\rfloor+1)^{2-1/r}}{\#\cB} q^{\frac{r+1}{4r^2}} (\log q)^{\frac{1}{2r}} T_i(q).    \end{align*}

Let $\phi^* = \phi(q)/q$, $B = B_i(r,q)$, $$\kappa = \left(\frac{(B-1)B^{1-\frac{1}{r}}}{2r(B+1)^{2-\frac{1}{r}}}\right)^{\frac{r}{r-1}},$$
and 
$$A = \frac{\kappa N}{B}.$$

Using \eqref{eq E sum bound}, Theorem \ref{weil}, Lemma \ref{v_A sum}, and Proposition \ref{real v_A square sum} yields
\footnote{In the first upper bound we actually have $(\lfloor B\rfloor+1)^{2-1/r}/\lfloor B\rfloor$; but the function $(x+1)^{2-1/r}/x$ is increasing for $x\ge2$, so we can replace $\lfloor B\rfloor$ by $B$ as shown.}
\begin{align}\label{eq together}
\begin{split}
    |S&_\chi(M, N)|
    \le
    \frac{1}{\#\cA\cdot \lfloor B\rfloor}\left(\#\cA N\right)^{1-\frac{1}{r}}
    \left(\#\cA^2 + 2AN\left(\frac{6}{\pi^2}\log{(A)}+\delta+\frac{\log{(A)}+1}{3A}\right)\right)^{\frac{1}{2r}}\times
    \\
    &\times\left(2r(4r)^{\omega(q)} B^{2r}m_r(q)\sqrt{q} + \frac{r^{2r}}{r!} B^r q\right)^\frac{1}{2r}+\frac{2C}{2-\frac{1}{r}} A^{1-\frac{1}{r}} \frac{(B+1)^{2-\frac{1}{r}}}{B} q^{\frac{r+1}{4r^2}} (\log q)^{\frac{1}{2r}} T_i(q)\\
    &\le \frac{N^{1-\frac{1}{r}}}{\#\cA^{\frac{1}{r}}\cdot (B-1)}\left(\left(\frac{\kappa N}{B}\phi^*+2^{\omega(q)-1}\right)^2 + \frac{2\kappa N^2}{B}\left(\frac{6}{\pi^2}\log{(A)}+\delta+\frac{\log{A}+1}{3A}\right)\right)^{\frac{1}{2r}}\times
    \\
    &\times\left(2r(4r)^{\omega(q)} B^{2r}m_r(q)\sqrt{q} + \frac{r^{2r}}{r!} B^r q\right)^\frac{1}{2r}+\frac{2C}{2-\frac{1}{r}}(\kappa N)^{1-\frac{1}{r}}\left(\frac{B+1}{B}\right)^{2-\frac{1}{r}}q^{\frac{r+1}{4r^2}} (\log q)^{\frac{1}{2r}} T_i(q).
\end{split}
\end{align}



Let
\begin{equation}\label{eq s}
   s=\frac{2C\kappa^{1-1/r}}{2-1/r}\left(\frac{B+1}{B}\right)^{2-1/r},
\end{equation}

\begin{equation*}
    P
    =
    \left(\frac{\kappa N}{B}\phi^*+2^{\omega(q)-1}\right)^2 + \frac{2\kappa N^2}{B}\left(\frac{6}{\pi^2}\log{(A)}+\delta+\frac{\log{A}+1}{3A}\right),
\end{equation*}
and
\begin{equation*}    
    Q
    =
    2r(4r)^{\omega(q)} B^{2r}m_r(q)\sqrt{q} + \frac{r^{2r}}{r!} B^r q.
\end{equation*}

We can write the upper bound for $|S_\chi(M, N)|$ in \eqref{eq together} as 
\begin{equation}
\label{eq together 2}
    \frac{N^{1-\frac{1}{r}}}{\#\cA^{\frac{1}{r}}\cdot (B-1)}P^{\frac{1}{2r}}\times
     Q^{\frac{1}{2r}}
    +
    N^{1-1/r}q^{\frac{r+1}{4r^2}}T_i(q)\left(\log q\right)^{\frac{1}{2r}}s.
\end{equation}
Factoring out $N^2/B$ from $P$, we deduce
\begin{align}
\label{eq 1/2r root P1}
    P
    &=
    \frac{N^2}{B}\left[ \frac{\left(\kappa\phi^*+\frac{B2^{\omega(q)-1}}{N} \right)^2}{B}
    +
    2\kappa\left(\frac{6}{\pi^2}\log{(A)}+\delta+\frac{\log{(A)}+1}{3A}\right)\right]
    \\
    \implies 
    P^{\frac{1}{2r}}
    &=
    N^{1/r}\left(\frac{\alpha}{B}\right)^{\frac{1}{2r}},\notag
\end{align}
where for simplicity, we let
\begin{equation}\label{eq alpha}
    \alpha
    =
    PB/N^2
    =
    \left[ \frac{\left(\kappa\phi^*+\frac{B2^{\omega(q)-1}}{N} \right)^2}{B}
    +
    2\kappa\left(\frac{6}{\pi^2}\log{(A)}+\delta+\frac{\log{(A)}+1}{3A}\right)\right].
\end{equation}

Let $$\beta = \frac{B}{r^2 q^{\frac{1}{2r}}}.$$
Then 
\begin{equation}
\label{eq 1/2r root P2}
    Q^{\frac{1}{2r}}=   q^{\frac{3}{4r}}r^{2}\beta^{1/2}\left(\frac{f_i}{r!}\right)^{\frac{1}{2r}},
\end{equation}
where
$f_1 = r$ and 
$$f_2 = (4r)^{\omega(q)}m_r(q)(r-1)+1.$$

Therefore, writing $B = r^2q^{\frac{1}{2r}}\beta$, we obtain from Equations \eqref{eq together 2}, \eqref{eq 1/2r root P1}, and \eqref{eq 1/2r root P2} that

\begin{equation}\label{helpful intermediate}
    \left| S_\chi(M, N) \right|
    \le
    \frac{N^{1-1/r}}{B-1}\left(\frac{N}{\#\cA}\right)^{\frac{1}{r}}q^{\frac{3}{4r}-\frac{1}{4r^2}}r^{2-\frac{1}{r}}\beta^{\frac{1}{2}-\frac{1}{2r}}\left(\frac{\alpha f_i}{r!}\right)^{\frac{1}{2r}}
        +N^{1-\frac{1}{r}}q^{\frac{r+1}{4r^2}}T_i(q)\log(q)^{\frac{1}{2r}}s.     
\end{equation}

Let
\begin{equation}\label{eq u}
u=\frac{\left(N/\#\cA\right)^{1/{r}}}{q^{\frac{1}{2r^2}}},
\end{equation}

and

\begin{equation}
\label{w}
    w
    = \frac{f_i}{r^2 \beta^{r+1} r! \log{q}}.
\end{equation}

By multiplying and dividing by $B$ and a few other manipulations we can go from \eqref{helpful intermediate} to 
\begin{equation}
\begin{split}
    \left| S_\chi(M, N) \right|
    &\le
    \left(\frac{B}{B-1}\right)N^{1-1/r}q^{(r+1)/{4r^2}}T_i(q)(\log q)^{\frac{1}{2r}}
    \left(
    \frac{u({\alpha w})^{\frac{1}{2r}}}{T_i(q)} 
    +s
    \right) \\
    &=
    \left(\frac{B}{B-1}\right)N^{1-1/r}q^{(r+1)/{4r^2}}(\log q)^{\frac{1}{2r}}
    \left(
    u(\alpha w)^{\frac{1}{2r}}
    +sT_i(q)
    \right),
\end{split} 
\end{equation}
as stated in the lemma.
\end{proof}

Using this lemma, we now prove our main theorems separately below.
\begin{proof}[Proof of Theorem \ref{main thm 1}]
Recall that we assume $q\ge \max\{10^{1143}, t_1(r)\}$ for $t_1(r)=e^{e^{a(r)}}$ with $a(r)$ given in the statement of Theorem \ref{main thm 1} (and Lemma \ref{B will be at least 10 if it kills us lemma}), and from
\eqref{eq T1}
$$T_1(q) = \left((4r)^{\omega(q)}m_r(q)\right)^{\frac{1}{2r}-\frac{1}{2r^2}}\left(\frac{q}{\phi(q)}\right)^{\frac{1}{r}}.$$
By Lemma \ref{lemma to set up the main theorems}
and comparing with \eqref{eq desired bound}, it suffices to show that
\begin{equation*}
u(\alpha w)^{\frac{1}{2r}} + sT_1(q) \le \frac{B-1}{B} C T_1(q).
\end{equation*}

Since $10^{1143}\ge (10r)^{18}$ for $r\le 10^{50}$ while $t_1(r) \ge (10r)^{18}$ for $r\ge 10^{50}$, we apply Lemma \ref{A lower bound lemma} and conclude that $A\phi^{*} \ge 2^{\omega(q)}$, hence
$$\frac{N}{\#\mathcal{A}}\le \frac{N}{A\phi^*-2^{\omega(q)-1}}=\frac{N}{A\phi^*}+\frac{2^{\omega(q)-1}N}{A\phi^*(A\phi^*-2^{\omega(q)-1})} \le \frac{2N}{A\phi^{*}}.$$
Therefore, \eqref{eq u} satisfies the bound
\begin{equation} \label{roughly u}
u=\frac{\left(N/\#\cA\right)^{1/{r}}}{q^\frac{1}{2r^2}} \le \frac{2^{1/r}}{q^{\frac{1}{2r^2}}}\left(\frac{N}{A\phi^*}\right)^{1/r}.
\end{equation}

From Lemma \ref{lem: alpha bound}, we have
\begin{align}\label{eq alpha upper bound}
\alpha \le \frac{32}{31}\kappa \log{q}.    
\end{align}
Using that $f_1 = r$ and \begin{equation}\label{beta B def}\beta = \left(\frac{r-1}{(2r)r!(4r)^{\omega(q)}m_r(q)}\right)^{\frac{1}{r}},
\end{equation}
we rewrite \eqref{w} as
\begin{equation}\label{bounding w}
    w = \frac{2^{1+1/r}(r\cdot r!)^{1/r}}{(r-1)^{1+1/r}\log q} \bigl(\left(4r\right)^{\omega(q)}m_r(q)\bigr)^{1+1/r}.
\end{equation}
Furthermore, using \eqref{beta B def}, $B = r^2q^{\frac{1}{2r}}\beta$, and $AB = \kappa N$, we have
\begin{equation}\label{bounding N over A}
\frac{N}{A q^{\frac{1}{2r}}}\left((4r)^{\omega(q)}m_r(q)\right)^{\frac{1}{r}} = \frac{r^2}{\kappa}\left(\frac{r-1}{r!(2r)}\right)^{\frac{1}{r}}.
\end{equation}
Therefore, combining \eqref{eq s}, \eqref{roughly u}, \eqref{eq alpha upper bound}, \eqref{bounding w}, and \eqref{bounding N over A},  we obtain
\begin{align*}
    &u(\alpha w)^{\frac{1}{2r}}+sT_1(q)\\ 
    &\le 2^{\frac{1}{r}}\left(\frac{N}{Aq^{\frac{1}{2r}}\phi^*}\right)^{\frac{1}{r}}\left(\frac{32}{31}\kappa \log{q}\right)^{\frac{1}{2r}}\left(\frac{2^{1+\frac{1}{r}}(r\cdot r!)^{\frac{1}{r}}}{(r-1)^{1+\frac{1}{r}}\log{q}}\right)^{\frac{1}{2r}}\left((4r)^{\omega(q)}m_r(q)\right)^{\frac{1}{2r}+\frac{1}{2r^2}}\\ &+\frac{2C\kappa^{1-\frac{1}{r}}}{2-\frac{1}{r}}\left(\frac{B+1}{B}\right)^{2-\frac{1}{r}}T_1(q)\\
    &= \left(2^{\frac{4}{r}-\frac{1}{2r^2}}31^{-\frac{1}{2r}}(r!)^{-\frac{1}{2r^2}}(r-1)^{\frac{1}{2r^2}-\frac{1}{2r}}r^{\frac{2}{r}-\frac{1}{2r^2}}\kappa^{-\frac{1}{2r}} +\frac{2C\kappa^{1-\frac{1}{r}}}{2-\frac{1}{r}}\left(\frac{B+1}{B}\right)^{2-\frac{1}{r}}\right)T_1(q).
\end{align*}
Therefore,
\begin{equation} \label{C lower bound}
C \ge \frac{2^{\frac{4}{r}-\frac{1}{2r^2}}31^{-\frac{1}{2r}}(r!)^{-\frac{1}{2r^2}}(r-1)^{\frac{1}{2r^2}-\frac{1}{2r}}r^{\frac{2}{r}-\frac{1}{2r^2}}\kappa^{-\frac{1}{2r}}}{\frac{B-1}{B}-\frac{2\kappa^{1-\frac{1}{r}}}{2-\frac{1}{r}}\left(\frac{B+1}{B}\right)^{2-\frac{1}{r}}}.
\end{equation}

Now to obtain the constants in Table \ref{Table 1}, first, by Lemmas \ref{B will be at least 10 if it kills us lemma} and \ref{B will be at least 10 if it kills us lemma 2},
\begin{equation}\label{B lower bound123}
B\ge 2^{2-\frac{2}{r}}r^2\left(\frac{1-\frac{1}{r}}{r!}\right)^{\frac{1}{r}}.
\end{equation}
Then, choosing $\kappa$ as in \eqref{choice of kappa}, we replace $B$ with $2^{2-2/r}r^2\left(\frac{1-\frac{1}{r}}{r!}\right)^{\frac{1}{r}}$ in $\kappa$ since the right side of \eqref{C lower bound} is decreasing as $B$ increases for $B\ge 2^{2-2/r}r^2\left(\frac{1-\frac{1}{r}}{r!}\right)^{\frac{1}{r}}$. Computing the resulting expression obtains the $C(r)$ column in Table \ref{Table 1}. 
For the $D(r)$ column in Table \ref{Table 1}, we notice that $q\rightarrow \infty$ implies $B\rightarrow \infty$, therefore $$\kappa\rightarrow \left(\frac{1}{2r}\right)^{\frac{r}{r-1}},$$
and we obtain
$$D(r) = \frac{2^{\frac{4}{r}-\frac{1}{2r^2}}31^{-\frac{1}{2r}}(r!)^{-\frac{1}{2r^2}}(r-1)^{\frac{1}{2r^2}-\frac{1}{2r}}r^{\frac{2}{r}-\frac{1}{2r^2}}(2r)^{\frac{1}{2r-2}}}{1-\frac{1}{2r-1}}.$$
\end{proof}

We prove Theorem \ref{main theorem 2} following the same procedure as the proof of Theorem \ref{main thm 1}.
\begin{proof}[Proof of Theorem \ref{main theorem 2}]
Recall that we assume $q\ge \max\{10^{1143}, t_2(r)\}$ for $t_2(r)=2^{4r-2}$, and from \eqref{eq T2}
$$T_2(q) = \left((4r)^{\omega(q)}m_r(q)\right)^{\frac{1}{2r}}\left(\frac{q}{\phi(q)}\right)^{\frac{1}{r}}.$$
By Lemma \ref{lemma to set up the main theorems}, it suffices to show
\begin{equation*}
u(\alpha w)^{\frac{1}{2r}} + sT_2(q) \le \frac{B-1}{B} C T_2(q).
\end{equation*}

Similar to the proof of Theorem \ref{main thm 1}, $10^{1143}\ge (10r)^{18}$ for $r\le 10^{50}$ while $t_2(r) \ge (10r)^{18}$ for $r\ge 10^{50}$, thus from \eqref{roughly u} and \eqref{eq alpha upper bound},
\begin{equation*}
    u \le \frac{2^{1/r}}{q^{\frac{1}{2r^2}}}\left(\frac{N}{A\phi^*}\right)^{1/r},
\end{equation*}
and 
\[
\alpha \le \frac{32}{31}\kappa \log{q}.
\]
Furthermore, using that $(4r)^{\omega(q)}\ge 8$,
\begin{equation}\label{f2}
f_2 = (4r)^{\omega(q)}m_r(q)(r-1)+1 \le \left(r-\frac{7}{8}\right)(4r)^{\omega(q)}m_r(q),
\end{equation}
and \begin{equation}\label{beta B def 2}\beta = \left(\frac{r-1}{(2r)r!}\right)^{\frac{1}{r}},
\end{equation}
therefore we bound $w$ in \eqref{w} as
\begin{equation}\label{bounding w 2}
    w \le \frac{2^{1+1/r}(r\cdot r!)^{1/r}}{(r-1)^{1+1/r}\log q} \left(\frac{r-\frac{7}{8}}{r}\right)\left(4r\right)^{\omega(q)}m_r(q).
\end{equation}
Using \eqref{beta B def 2}, $B = r^2q^{\frac{1}{2r}} \beta$, and $AB = \kappa N$,  we get
\begin{equation}\label{bounding N over A 2}
\frac{N}{A q^{\frac{1}{2r}}} = \frac{r^2}{\kappa}\left(\frac{r-1}{r!(2r)}\right)^{\frac{1}{r}}.
\end{equation}
Therefore, combining \eqref{eq s}, \eqref{roughly u}, \eqref{eq alpha upper bound}, \eqref{bounding w 2}, and \eqref{bounding N over A 2}, we get 
\begin{align*}
    &u(\alpha w)^{\frac{1}{2r}}+sT_2(q)\\ 
    &\le 2^{\frac{1}{r}}\left(\frac{N}{Aq^{\frac{1}{2r}}\phi^*}\right)^{\frac{1}{r}}\left(\frac{32}{31}\kappa \log{q}\right)^{\frac{1}{2r}}\left(\frac{2^{1+\frac{1}{r}}(r\cdot r!)^{\frac{1}{r}}}{(r-1)^{1+\frac{1}{r}}\log{q}}\left(\frac{r-\frac{7}{8}}{r}\right)\right)^{\frac{1}{2r}}\left((4r)^{\omega(q)}m_r(q)\right)^{\frac{1}{2r}}\\ &+\frac{2C\kappa^{1-\frac{1}{r}}}{2-\frac{1}{r}}\left(\frac{B+1}{B}\right)^{2-\frac{1}{r}}T_2(q)\\
    &= \left(2^{\frac{4}{r}-\frac{1}{2r^2}}31^{-\frac{1}{2r}}(r!)^{-\frac{1}{2r^2}}(r-1)^{\frac{1-r}{2r^2}}r^{\frac{3r-1}{2r^2}}\left(r-\frac{7}{8}\right)^{\frac{1}{2r}}\kappa^{-\frac{1}{2r}} +\frac{2C\kappa^{1-\frac{1}{r}}}{2-\frac{1}{r}}\left(\frac{B+1}{B}\right)^{2-\frac{1}{r}}\right)T_2(q).
\end{align*}
Therefore,
\begin{equation} \label{C* lower bound}
C \ge \frac{2^{\frac{4}{r}-\frac{1}{2r^2}}31^{-\frac{1}{2r}}(r!)^{-\frac{1}{2r^2}}(r-1)^{\frac{1}{2r^2}-\frac{1}{2r}}r^{\frac{3}{2r}-\frac{1}{2r^2}}\left(r-\frac{7}{8}\right)^{\frac{1}{2r}}\kappa^{-\frac{1}{2r}}}{\frac{B-1}{B}-\frac{2\kappa^{1-\frac{1}{r}}}{2-\frac{1}{r}}\left(\frac{B+1}{B}\right)^{2-\frac{1}{r}}}.
\end{equation}
We make the same choices for $\kappa$ as in the proof of Theorem \ref{main thm 1} which leads to the values of $C^*(r)$ in Table \ref{Table 1}.
\end{proof}


\section*{Acknowledgments}

The authors thank the PIMS CRG ``$L$-functions in Analytic Number Theory'' and the BIRS facilities at the University of British Columbia--Okanagan, who together ran the Inclusive Paths in Explicit Number Theory Summer School where this project started. The third author was supported in part by a Natural Sciences and Engineering Council of Canada Discovery Grant. The fourth and fifth authors also thank FONDOCyT, grant number 2022-1D1-085, for supporting this research. We also wish to thank the anonymous referee, specially for suggesting Lemma \ref{phi sum lemma} and other observations that led to some improvements in the constants in our main theorems.

\bibliographystyle{amsplain}
\bibliography{bsample}

@article {Elkies,
    AUTHOR = {Elkies, Noam D.},
     TITLE = {Rank of an elliptic curve and 3-rank of a quadratic field via
              the {B}urgess bounds},
   JOURNAL = {Res. Number Theory},
  FJOURNAL = {Research in Number Theory},
    VOLUME = {11},
      YEAR = {2025},
    NUMBER = {3},
     PAGES = {Paper No. 70, 14},
      ISSN = {2522-0160,2363-9555},
   MRCLASS = {11G05 (11L40 11R29 11Y40 11Y50)},
  MRNUMBER = {4933368},
       DOI = {10.1007/s40993-024-00601-x},
       URL = {https://doi.org/10.1007/s40993-024-00601-x},
}

@article {JRT-2023,
    AUTHOR = {Johnston, D. R. and Ramar\'e, O. and Trudgian, T.},
     TITLE = {An explicit upper bound for {$L(1,\chi)$} when {$\chi$} is
              quadratic},
   JOURNAL = {Res. Number Theory},
  FJOURNAL = {Research in Number Theory},
    VOLUME = {9},
      YEAR = {2023},
    NUMBER = {4},
     PAGES = {Paper No. 72, 20},
      ISSN = {2522-0160,2363-9555},
   MRCLASS = {11M20},
  MRNUMBER = {4649452},
MRREVIEWER = {Peter\ Shiu},
       DOI = {10.1007/s40993-023-00476-4},
       URL = {https://doi.org/10.1007/s40993-023-00476-4},
}

@article {COT-2016,
    AUTHOR = {Cohen, Stephen D. and Oliveira e Silva, Tom\'as and Trudgian,
              Tim},
     TITLE = {On {G}rosswald's conjecture on primitive roots},
   JOURNAL = {Acta Arith.},
  FJOURNAL = {Acta Arithmetica},
    VOLUME = {172},
      YEAR = {2016},
    NUMBER = {3},
     PAGES = {263--270},
      ISSN = {0065-1036,1730-6264},
   MRCLASS = {11L40 (11A07)},
  MRNUMBER = {3460815},
MRREVIEWER = {Ke\ Gong},
       DOI = {10.4064/aa8109-12-2015},
       URL = {https://doi.org/10.4064/aa8109-12-2015},
}

@article {burgess86,
    AUTHOR = {Burgess, D. A.},
     TITLE = {The character sum estimate with {$r=3$}},
   JOURNAL = {J. London Math. Soc. (2)},
  FJOURNAL = {Journal of the London Mathematical Society. Second Series},
    VOLUME = {33},
      YEAR = {1986},
    NUMBER = {2},
     PAGES = {219--226},
      ISSN = {0024-6107,1469-7750},
   MRCLASS = {11L40},
  MRNUMBER = {838632},
MRREVIEWER = {E.\ Stankus},
       DOI = {10.1112/jlms/s2-33.2.219},
       URL = {https://doi.org/10.1112/jlms/s2-33.2.219},
}

@article {francis2021,
    AUTHOR = {Francis, Forrest J.},
     TITLE = {An investigation into explicit versions of {B}urgess' bound},
   JOURNAL = {J. Number Theory},
  FJOURNAL = {Journal of Number Theory},
    VOLUME = {228},
      YEAR = {2021},
     PAGES = {87--107},
      ISSN = {0022-314X,1096-1658},
   MRCLASS = {11L40 (11Y60)},
  MRNUMBER = {4271811},
       DOI = {10.1016/j.jnt.2021.03.018},
       URL = {https://doi.org/10.1016/j.jnt.2021.03.018},
}

@article {mcgown2012,
    AUTHOR = {McGown, Kevin J.},
     TITLE = {Norm-{E}uclidean cyclic fields of prime degree},
   JOURNAL = {Int. J. Number Theory},
  FJOURNAL = {International Journal of Number Theory},
    VOLUME = {8},
      YEAR = {2012},
    NUMBER = {1},
     PAGES = {227--254},
      ISSN = {1793-0421,1793-7310},
   MRCLASS = {11R16 (11L40 11R32 11R80 11Y40)},
  MRNUMBER = {2887892},
MRREVIEWER = {J\"urgen\ Ritter},
       DOI = {10.1142/S1793042112500133},
       URL = {https://doi.org/10.1142/S1793042112500133},
}

@article {JKL,
    AUTHOR = {Jain-Sharma, Niraek and Khale, Tanmay and Liu, Mengzhen},
     TITLE = {Explicit {B}urgess bound for composite moduli},
   JOURNAL = {Int. J. Number Theory},
  FJOURNAL = {International Journal of Number Theory},
    VOLUME = {17},
      YEAR = {2021},
    NUMBER = {10},
     PAGES = {2207--2219},
      ISSN = {1793-0421,1793-7310},
   MRCLASS = {11L40},
  MRNUMBER = {4322829},
MRREVIEWER = {Matteo\ Bordignon},
       DOI = {10.1142/S1793042121500834},
       URL = {https://doi.org/10.1142/S1793042121500834},
}

@article {FroSo,
    AUTHOR = {Frolenkov, D. A. and Soundararajan, K.},
     TITLE = {A generalization of the {P}\'{o}lya-{V}inogradov inequality},
   JOURNAL = {Ramanujan J.},
  FJOURNAL = {Ramanujan Journal. An International Journal Devoted to the
              Areas of Mathematics Influenced by Ramanujan},
    VOLUME = {31},
      YEAR = {2013},
    NUMBER = {3},
     PAGES = {271--279},
      ISSN = {1382-4090},
   MRCLASS = {11L40 (11A25 11L03 11L07)},
  MRNUMBER = {3081668},
MRREVIEWER = {Shuguang Li},
       DOI = {10.1007/s11139-012-9462-y},
       URL = {https://doi.org/10.1007/s11139-012-9462-y},
}

@article {RS1962,
    AUTHOR = {Rosser, J. Barkley and Schoenfeld, Lowell},
     TITLE = {Approximate formulas for some functions of prime numbers},
   JOURNAL = {Illinois J. Math.},
  FJOURNAL = {Illinois Journal of Mathematics},
    VOLUME = {6},
      YEAR = {1962},
     PAGES = {64--94},
      ISSN = {0019-2082},
   MRCLASS = {10.42},
  MRNUMBER = {0137689 (25 \#1139)},
MRREVIEWER = {B. K. Ghosh},
}

@book {IK2004,
    AUTHOR = {Iwaniec, Henryk and Kowalski, Emmanuel},
     TITLE = {Analytic number theory},
    SERIES = {American Mathematical Society Colloquium Publications},
    VOLUME = {53},
 PUBLISHER = {American Mathematical Society},
   ADDRESS = {Providence, RI},
      YEAR = {2004},
     PAGES = {xii+615},
      ISBN = {0-8218-3633-1},
   MRCLASS = {11-02 (11Fxx 11Lxx 11Mxx 11Nxx)},
  MRNUMBER = {2061214 (2005h:11005)},
MRREVIEWER = {K. Soundararajan},
}

@article {Bur1962,
    AUTHOR = {Burgess, D. A.},
     TITLE = {On character sums and primitive roots},
   JOURNAL = {Proc. London Math. Soc. (3)},
  FJOURNAL = {Proceedings of the London Mathematical Society. Third Series},
    VOLUME = {12},
      YEAR = {1962},
     PAGES = {179--192},
      ISSN = {0024-6115},
   MRCLASS = {10.41},
  MRNUMBER = {0132732 (24 \#A2569)},
MRREVIEWER = {L. Carlitz},
}

@article {Bur1963,
    AUTHOR = {Burgess, D. A.},
     TITLE = {On character sums and {$L$}-series. {II}},
   JOURNAL = {Proc. London Math. Soc. (3)},
  FJOURNAL = {Proceedings of the London Mathematical Society. Third Series},
    VOLUME = {13},
      YEAR = {1963},
     PAGES = {524--536},
      ISSN = {0024-6115},
   MRCLASS = {10.41},
  MRNUMBER = {0148626 (26 \#6133)},
MRREVIEWER = {L. Carlitz},
}

@article {Bur19632,
    AUTHOR = {Burgess, D. A.},
     TITLE = {A note on the distribution of residues and non-residues},
   JOURNAL = {J. London Math. Soc.},
  FJOURNAL = {Journal of the London Mathematical Society. Second Series},
    VOLUME = {38},
      YEAR = {1963},
     PAGES = {253--256},
      ISSN = {0024-6107},
   MRCLASS = {10.41 (10.43)},
  MRNUMBER = {0148628 (26 \#6135)},
MRREVIEWER = {D. J. Lewis},
}

@article {Boo2006,
    AUTHOR = {Booker, Andrew R.},
     TITLE = {Quadratic class numbers and character sums},
   JOURNAL = {Math. Comp.},
  FJOURNAL = {Mathematics of Computation},
    VOLUME = {75},
      YEAR = {2006},
    NUMBER = {255},
     PAGES = {1481--1492 (electronic)},
      ISSN = {0025-5718},
     CODEN = {MCMPAF},
   MRCLASS = {11R47 (11L40 11M20 11R11 11R29 11Y35)},
  MRNUMBER = {2219039 (2008a:11140)},
MRREVIEWER = {John B. Friedlander},
       DOI = {10.1090/S0025-5718-06-01850-3},
       URL = {http://dx.doi.org/10.1090/S0025-5718-06-01850-3},
}

@article {Weil,
    AUTHOR = {Weil, Andr\'{e}},
     TITLE = {On some exponential sums},
   JOURNAL = {Proc. Nat. Acad. Sci. U. S. A.},
  FJOURNAL = {Proceedings of the National Academy of Sciences of the United
              States of America},
    VOLUME = {34},
      YEAR = {1948},
     PAGES = {204--207},
      ISSN = {0027-8424},
   MRCLASS = {10.0X},
  MRNUMBER = {0027006},
MRREVIEWER = {H. D. Kloosterman},
       DOI = {10.1073/pnas.34.5.204},
       URL = {https://doi.org/10.1073/pnas.34.5.204},
}

@article {ET1,
    AUTHOR = {Trevi\~{n}o, Enrique},
     TITLE = {The least {$k$}-th power non-residue},
   JOURNAL = {J. Number Theory},
  FJOURNAL = {Journal of Number Theory},
    VOLUME = {149},
      YEAR = {2015},
     PAGES = {201--224},
      ISSN = {0022-314X},
   MRCLASS = {11E12 (11D09 11E20)},
  MRNUMBER = {3296008},
MRREVIEWER = {Ronald J. Evans},
       DOI = {10.1016/j.jnt.2014.10.019},
       URL = {https://doi.org/10.1016/j.jnt.2014.10.019},
}

@article {ET2,
    AUTHOR = {Trevi\~{n}o, Enrique},
     TITLE = {The {B}urgess inequality and the least {$k$}th power
              non-residue},
   JOURNAL = {Int. J. Number Theory},
  FJOURNAL = {International Journal of Number Theory},
    VOLUME = {11},
      YEAR = {2015},
    NUMBER = {5},
     PAGES = {1653--1678},
      ISSN = {1793-0421},
   MRCLASS = {11L40 (11A15 11Y60)},
  MRNUMBER = {3376232},
MRREVIEWER = {John H. Loxton},
       DOI = {10.1142/S1793042115400163},
       URL = {https://doi.org/10.1142/S1793042115400163},
}

@article {burgess:1962,
    AUTHOR = {Burgess, D. A.},
     TITLE = {On character sums and {$L$}-series},
   JOURNAL = {Proc. London Math. Soc. (3)},
  FJOURNAL = {Proceedings of the London Mathematical Society. Third Series},
    VOLUME = {12},
      YEAR = {1962},
     PAGES = {193--206},
      ISSN = {0024-6115,1460-244X},
   MRCLASS = {10.41},
  MRNUMBER = {132733},
MRREVIEWER = {L.\ Carlitz},
       DOI = {10.1112/plms/s3-12.1.193},
       URL = {https://doi.org/10.1112/plms/s3-12.1.193},
}

@Article{Nicolas-Robin1983,
 Author = {Nicolas, J.-L. and Robin, G.},
 Title = {Explicit upper estimates for the number of divisors of {{\(n\)}}},
 FJournal = {Canadian Mathematical Bulletin},
 Journal = {Can. Math. Bull.},
 ISSN = {0008-4395},
 Volume = {26},
 Pages = {485--492},
 Year = {1983},
 Language = {French},
 DOI = {10.4153/CMB-1983-078-5},
 zbMATH = {3783121},
 Zbl = {0497.10034}
}

@article{BrokerSutherland2010,
  author  = {Reinier Br{\"o}ker and Andrew V. Sutherland},
  title   = {An Explicit Height Bound for the Classical Modular Polynomial},
  journal = {The Ramanujan Journal},
  volume  = {22},
  number  = {3},
  pages   = {293--313},
  year    = {2010},
  doi     = {10.1007/s11139-010-9231-8}
}

\end{document}